\documentclass[12pt]{amsart}
\usepackage{geometry}
\usepackage{amsmath,amsthm,amsfonts,amssymb,bm,graphicx, mathrsfs, mathabx }
\usepackage{verbatim}

\usepackage
{hyperref}
\hypersetup{colorlinks=true,citecolor=blue,linkcolor=blue,urlcolor=blue}

\theoremstyle{plain}
\newtheorem{theorem}{Theorem}

\newtheorem{prop}[theorem]{Proposition}
\numberwithin{equation}{section}

\theoremstyle{definition}

\usepackage{subcaption}
\usepackage{graphicx}

\makeatletter
\DeclareRobustCommand\widecheck[1]{{\mathpalette\@widecheck{#1}}}
\def\@widecheck#1#2{%
    \setbox\z@\hbox{\m@th$#1#2$}%
    \setbox\tw@\hbox{\m@th$#1%
       \widehat{%
          \vrule\@width\z@\@height\ht\z@
          \vrule\@height\z@\@width\wd\z@}$}%
    \dp\tw@-\ht\z@
    \@tempdima\ht\z@ \advance\@tempdima2\ht\tw@ \divide\@tempdima\thr@@
    \setbox\tw@\hbox{%
       \raise\@tempdima\hbox{\scalebox{1}[-1]{\lower\@tempdima\box
\tw@}}}%
    {\ooalign{\box\tw@ \cr \box\z@}}}
\makeatother

\begin{document}

\author{Valentin Blomer}
\address{University of Bonn, Mathematisches Institut, Endenicher Allee 60, D-53115 Bonn, Germany} \email{blomer@math.uni-bonn.de}

\author{Maksym Radziwi\l\l}
\address{
  Caltech,
  Department of Mathematics,
  1200 E California Blvd,
  Pasadena, CA, 91125}
\email{maksym.radziwill@gmail.com}

\title{Small gaps in the spectrum of tori: asymptotic formulae}
 

\thanks{The first author was partially supported by a SNF-DFG lead agency grant BL 915/2-2 and BL 915/5-1.
The second author acknowledges support of a Sloan Fellowship and NSF grant DMS-1902063}

\begin{abstract} We establish  an asymptotic formula, uniformly down to the Planck scale, for the number of small gaps between the first $N$ eigenvalues of the Laplacian on almost all flat tori and also on almost all rectangular flat tori. \end{abstract}


\setcounter{tocdepth}{2}  

\maketitle 

\section{Introduction}

What is the distribution of values at integer arguments of a randomly chosen positive binary quadratic form? This number theoretic question has the following well-known dynamical interpretation.   If $\Lambda\subseteq \Bbb{R}^2$ is a lattice of rank 2, then the numbers $4\pi^2 \| \bm \omega \|^2$, $\bm \omega \in \Lambda^{\ast}$, are the eigenvalues of the Laplacian on $\Bbb{R}^2/\Lambda$ which are given by the values of a positive binary quadratic form at integer arguments. The equations of motions of a free particle moving through $\Lambda$ are integrable, thus according to the conjectures of Berry-Tabor the energy levels of the corresponding quantized system (i.e.,\ $4\pi \| \bm \omega \|^2$, $\bm \omega \in \Lambda^{\ast}$) should behave like a sequence of points coming from a Poisson process, at least for generic $\Lambda$.


The first to investigate this phenomenon systematically in the case of quadratic forms was Sarnak   \cite{Sa} who showed that for almost all (in a Lebesgue sense) quadratic forms the pair correlation is Poissonian. To state this more precisely, let us fix some notation.  Given $\bm \alpha = (\alpha_1, \alpha_2, \alpha_3) \in \Bbb{R}^3$ with $4 \alpha_1 \alpha_3 > \alpha_2^2$ and $\alpha_1 > 0$, let $
q_{\bm \alpha}(m,n) = \alpha_1 m^2 + \alpha_2 m n + \alpha_3 n^2$ denote the corresponding positive binary quadratic form. Without loss of generality we may assume that $q_{\bm \alpha}$ is reduced, so that we can restrict ${\bm \alpha } \in \mathcal{D}$ with 
$$\mathcal{D} = \{(\alpha_1, \alpha_2, \alpha_3) \in \Bbb{R}^3 \mid 0 \leq \alpha_2 \leq \alpha_1 \leq \alpha_3\}. $$
 The precise measure that we choose on $\mathcal{D}$ is relatively unimportant, but it is  most natural to choose the ${\rm GL}_2(\Bbb{R})$-invariant hyperbolic measure
$$d_{\text{hyp}}{\bm \alpha} = \frac{d\bm \alpha}{\pi^3 D(\bm \alpha)^3} = \frac{d\alpha_1 \, d\alpha_2\, d\alpha_3}{(4\alpha_1\alpha_3 - \alpha^2)^{3/2}}.$$

Each quadratic form $q_{\bm \alpha}$ has the automorphism $(m, n) \mapsto (-m, -n)$, so each positive eigenvalue occurs with multiplicity at least 2. It is therefore natural to desymmetrize the spectrum and consider only the values $q_{\bm \alpha}(m, n)$ with $m > 0$ or $m=0$ and $n \geq 0$.  We denote by $0 < \Lambda_1 \leq \Lambda_2 \leq \ldots $ the ordered set of values 
\begin{equation*}
\frac{q_{\bm \alpha}(m,n) }{D(\bm \alpha)}   , \quad D(\bm \alpha) = \frac{1}{\pi} \sqrt{4 \alpha_1 \alpha_3 - \alpha_2^2},
\end{equation*}
where $m > 0$ or $m = 0$ and $n \geq 0$. We usually suppress the dependence on ${\bm \alpha}$ of the numbers $\Lambda_j$. 
Asymptotically the average spacing between the $\Lambda_i$ is one, thus we think of $\Lambda_i$ as the properly rescaled multi-set of  eigenvalues of the  Laplacian on a suitable lattice. For an interval $I \subseteq \Bbb{R}$ and $N \geq 1$ we write
\begin{equation}\label{defP}
P({\bm \alpha}, N, I) = \frac{1}{N }\#\big\{(j, k) \mid \Lambda_j, \Lambda_k \leq N, j \not = k,   \Lambda_j - \Lambda_k\in I\big\}
\end{equation}
and denote by $\mu(I)$ the length of $I$. Sarnak \cite[Theorem 1]{Sa} shows that almost all $\bm \alpha$ satisfy
$$P({\bm \alpha}, N, I)  \sim \mu(I)$$
for any fixed interval $I$ as $N \rightarrow \infty$. That this holds even for all \emph{diophantine} ${\bm \alpha}$ was proved in \cite[Theorem 1.7]{EMM}. It is clear that it cannot hold for all ${\bm \alpha}$, for instance it is clearly wrong for all integral forms $q_{\bm \alpha}$ and any nonempty $I$ not containing zero of size strictly less than $D(\bm \alpha)^{-1}$. 

Sarnak's result is an effective asymptotic formula that comes with a(n unspecified) power saving in the error term, so it is clear that one can shrink the interval $I$ a little bit with $N$.  Our first result shows that the asymptotic formula remains true for almost all $\bm \alpha$ even for $\mu(I)$ as small as $N^{-1+\varepsilon}$. This is, up to the value of $\varepsilon$, the smallest scale at which we expect that gaps exist at all, it corresponds to the Planck scale. In this sense Theorem \ref{thm1}  is best possible. 

\begin{theorem}\label{thm1}
Let $\eta > 0$ be given.  There exists a subset $\mathcal{E} \subseteq \mathcal{D}$ of measure zero with the following property: for all ${\bm \alpha} \in \mathcal{D} \setminus \mathcal{E}$ we have
\begin{equation}\label{asymp}
 P({\bm \alpha}, N, [0, \Delta]) = (1 + o(1)) \Delta,  \quad N \rightarrow \infty, 
 \end{equation}
uniformly in 
\begin{equation}\label{NDelta}
N^{-1+\eta} \leq \Delta  \leq N^{-\eta}.
\end{equation}
\end{theorem}

The restriction $\Delta \leq N^{-\eta}$ could easily be removed. We included it for convenience to streamline the argument as the main interest is certainly the case of small $\Delta$. \\

Our next result concerns the thin subset of rectangular tori where $\alpha_2  = 0$, i.e.\ the value distribution of \emph{diagonal} quadratic forms. It was shown in \cite[Theorem 1.2]{BBRR} that almost all rectangular tori have pairs of eigenvalues of size at most $N$ with difference at most $N^{-1+\varepsilon}$. Here we upgrade the mere existence of small gaps to an asymptotic formula for its cardinality. 

Diagonal forms have four symmetries generated by $(m, n) \mapsto (m, -n)$ and $(m, n) \mapsto (-m, n)$, and we denote by $0 < \Lambda_1 \leq \Lambda_2 \leq \ldots$ the ordered sequence of values of
$$\frac{\pi}{4\sqrt{\alpha_1\alpha_3}}   q_{\bm \alpha}(m, n), \quad m > 0, n \geq 0$$
and accordingly define $P({\bm \alpha}, N, I)$ as in \eqref{defP}. 
Let  $\mathcal{R}$ be the set of $ (\alpha_1, \alpha_3) \in \Bbb{R}^2_{>0}$ which is naturally equipped with the measure $\frac{d\alpha_1\, d\alpha_3}{\alpha_1\alpha_3}$.

\begin{theorem}\label{thm2} Let $\eta > 0$ be given.  There exists a subset $\mathcal{F} \subseteq \mathcal{R}$ of measure zero with the following property: for all ${\bm \alpha} \in \mathcal{R} \setminus \mathcal{F}$ we have
\begin{equation*}
 P({\bm \alpha}, N, [0, \Delta]) = (1 + o(1)) \Delta,  \quad N \rightarrow \infty, 
 \end{equation*}
uniformly in 
\begin{equation}\label{NDelta1}
N^{-1+\eta} \leq \Delta  \ll 1.
\end{equation}
\end{theorem}


The proofs of Theorems \ref{thm1} and \ref{thm2} use a variety of techniques. Both of them start with Fourier analysis and transform the problem at hand into a diophantine question. The arithmetic part of the proof of Theorem \ref{thm1} is mainly based on lattice point arguments and the geometry of numbers. Theorem \ref{thm2} uses more advanced machinery. The desired asymptotic formula follows without much difficulty from the Generalized Riemann Hypothesis, or even from a Lindel\"of-type bound for the 8th moment of the Riemann zeta function on the half-line. In \cite{BBRR} the use of GRH was avoided by introducing an additional bilinear structure (and hence a second set of variables) along with the  best known bounds for the Riemann zeta function close to the one-line. This comes at the price of losing density in the asymptotic formula and returns only a lower bound for $P({\bm \alpha}, N, [0, \Delta])$. Therefore we need a new idea. The plan is to restrict the second set of variables to primes and employ ideas from \cite{MR} together with an analysis of numbers without small and large prime factors.

The problem of determining the smallest $\Delta$ for which
$$\liminf_{N \rightarrow \infty} P({\bm \alpha}, N, [0, \Delta]) > 0 \qquad \text{almost surely in } \bm \alpha$$ has attracted some attention, in the context of Theorem \ref{thm2}. Recently Aistleitner, El-Baz and Munsch \cite{AEM} showed that for almost $\bm \alpha$ there exist gaps that are at most $ (\log N)^{2c} / N$ with $c = 1 - \frac{1 + \log\log 2}{\log 2} \approx 0.086$ the Erd\H{o}s-Tenenbaum-Ford constant. In this direction we note that with more effort Theorem \ref{thm2} can be shown to still hold for $\Delta = (\log N)^{A} / N$ and $A > 0$ some large fixed constant.

\section{Proof of Theorem \ref{thm1}}

Since $\mathcal{D}$ can be covered by countably many compact sets, it suffices to consider a compact subset $\mathcal{D}_0 \subseteq \mathcal{D}$ and show that almost all ${\bm \alpha} \in  \mathcal{D}_0 $ satisfy \eqref{asymp}. Next we observe that $P({\bm \alpha}, N, I) = P(\lambda{\bm \alpha}, N, I)$ for every $\lambda > 0$, so we can de-homogenize by setting  $\alpha_2 = 1$, as the set of ${\bm \alpha} \in \mathcal{D}_0$ with $\alpha_2 = 0$  has measure 0.  From   now on we write  
$$q_{\bm \alpha}(m,n) = \alpha_1 m^2 +   m n + \alpha_3 n^2$$
where ${\bm \alpha} = (\alpha_1, 1, \alpha_3)$ is contained in a compact domain $\mathcal{D}^{\ast}$ inside 
\begin{equation}\label{inside}
 \{(\alpha_1, 1, \alpha_3) \in \Bbb{R}^2 \mid 1 \leq \alpha_1 \leq \alpha_3\}. 
 \end{equation}
Correspondingly we write
$$d_{\text{hyp}}^{\ast} {\bm \alpha} = \frac{d\alpha_1 \, d\alpha_3}{(4\alpha_1 \alpha_3 - 1)^{3/2}} = \frac{d{\bm \alpha}}{\pi^3 D(\bm \alpha)}$$
where now $D({\bm \alpha}) = D((\alpha_1, 1, \alpha_3)) = \pi^{-1} (4\alpha_1\alpha_3 - 1)^{1/2}.$ 
Let $V, W$ be  fixed smooth, real-valued functions with compact support in  $(-\infty, \infty)$ respectively. 
For $M, T \geq 1$ we define
$$ \mathcal{G}_{\bm \alpha}(M, T) :=  \frac{1}{4}\sum_{\substack{x_1, x_2, y_1, y_2 \in \Bbb{Z}\\(x_1, y_1) \not= \pm (x_2, y_2)}} W\Big( T\cdot \frac{q_{\bm \alpha}(x_1, y_1) - q_{\bm \alpha}(x_2, y_2)}{D({\bm \alpha}) }\Big) V\Big( \frac{q_{\bm \alpha}(x_1, y_1)}{M^2 D({\bm \alpha})}\Big).$$
As the argument of $V$ is non-negative, we define
\begin{equation}\label{vast}
V^{\ast}(x) = \delta_{x \geq 0} V(x).
\end{equation}
As a precursor to  Theorem \ref{thm1} we show the following proposition.
\begin{prop}\label{prop3} Fix $\eta > 0$ and suppose that
\begin{equation}\label{MT}
 M^{\eta} \leq T \leq M^{2-\eta}.
\end{equation}
Then for all $\varepsilon > 0$ we have 
$$\int_{\mathcal{D}^{\ast}}  \Big|\mathcal{G}_{\bm \alpha}(M, T) - \widehat{V^{\ast}}(0) \widehat{W}(0) \frac{M^2}{T} \Big|^2 d^{\ast}_{\text{{\rm hyp}}}{\bm \alpha}\ll_{\varepsilon} \frac{M^{4 - \frac{1}{44}+\varepsilon}}{T^2} + \frac{M^{2+\varepsilon}}{T} .$$
\end{prop}
Here and in the following we denote by $\widehat{f}$ the Fourier transform of $f$.
 We postpone the proof of Proposition \ref{prop3} to   Section \ref{proof3} and complete the proof of Theorem \ref{thm1}. For $0 < \delta < 1$ and $M, T$ as in \eqref{MT} let 
 $\mathcal{S}^+_{\delta}(M, T^{-1})$ be the set of ${\bm \alpha}$ such that
 $$ \#\big\{(j, k) \mid \Lambda_j, \Lambda_k \leq M^2, j \not = k, 0 \leq   \Lambda_j - \Lambda_k \leq T^{-1}\big\} \geq (1 + \delta)  \frac{M^2}{T} . $$
 We specialize $V, W$ to be smooth, non-negative functions such that $V, W \geq 1$ on $[0, 1]$ and $1 \leq \widehat{V}^{\ast}(0), \widehat{W}(0) \leq 1 + \delta/3$. Then $\mathcal{S}^+_{\delta}(M, T)$ is contained in the set of ${\bm \alpha}$ such that
 $$\mathcal{G}_{\bm \alpha}(M, T) \geq (1 + \delta) \frac{M^2}{T},$$
 so that $$\mu_{\text{hyp}}(\mathcal{S}^+_{\delta}(M, T^{-1})) \ll  \delta^{-2} M^{-\eta/2} $$
 if \eqref{MT} holds and $\eta \leq 1/23$. In the same way we can bound the measure of the set  $\mathcal{S}^-_{\delta}(M, T^{-1})$  of ${\bm \alpha}$ such that
 $$ \#\big\{(j, k) \mid \Lambda_j, \Lambda_k \leq M^2, j \not = k, 0 \leq   \Lambda_j - \Lambda_k \leq T^{-1}\big\} \leq (1 - \delta)  \frac{M^2}{T} . $$
 Hence if  $\mathcal{S}_{\delta}(N, \Delta)$ is the set of ${\bm \alpha}$ such that
 $$\Big| \#\big\{(j, k) \mid \Lambda_j, \Lambda_k \leq N, j \not = k, 0 \leq   \Lambda_j - \Lambda_k \leq \Delta\big\}  - N\Delta\Big| \geq \delta N \Delta,  $$
we conclude that 
 $$\mu_{\text{hyp}}(\mathcal{S}_{\delta}(N, \Delta)) \ll  \delta^{-2} M^{-\eta/2} $$
 uniformly in the region \eqref{NDelta}. Now let $\mathcal{S}_{\delta}$ be the set of ${\bm \alpha}$ such that there exists a pair of sequences $N_j, \Delta_j$ with $N_j \rightarrow \infty$ and  $(N_j, \Delta_j)$ satisfying \eqref{NDelta} such that ${\bm \alpha} \in \mathcal{S}_{\delta}(N_j, \Delta_j)$ for all $j$. For approximation purposes we now consider the special sequences $N^{\ast}_m= (1+\delta^2)^m$, and $\Delta^{\ast}_n = (1+\delta^2)^{-n}$, where $m, n \in \Bbb{N}$.  If $\delta$ is sufficiently small, then for each $j$ there exists a pair $m, n$ such $|N_j - N^{\ast}_m| \ll \delta^2 N_j $ and $|\Delta_j - \Delta^{\ast}_n| \ll \delta^2 \Delta_j $ and so 
 $\mathcal{S}_{\delta}(N_j, \Delta_j) \subseteq\mathcal{S}_{\delta/2}(N^{\ast}_m, \Delta^{\ast}_n)$. Here we have necessarily $n \ll \log N^{\ast}_m$. We conclude that 
 $$\mathcal{S}_{\delta}(N_j, \Delta_j) \subseteq \bigcup_{n \ll \log N^{\ast}_m} \mathcal{S}_{\delta/2}(N^{\ast}_m, \Delta^{\ast}_n) =: \mathcal{S}^{\ast}_{\delta/2}(N^{\ast}_m),$$
 and hence $\mathcal{S}_{\delta} \subseteq \limsup_m \mathcal{S}^{\ast}_{\delta/2}(N^{\ast}_m)$. As 
 $$\mu_{\text{hyp}}(\mathcal{S}^{\ast}_{\delta/2}(N^{\ast}_m)) \ll (N^{\ast}_m)^{-\eta/2}\log N^{\ast}_m = (1 + \delta^2)^{-\eta m/2} \log ( 1+ \delta^2)^m$$
 we conclude from the Borel-Cantelli lemma that $\mu_{\text{hyp}}(\mathcal{S}_{\delta}) = 0$. Therefore, the set
 $$\bigcup_{n} \mathcal{S}_{1/n}$$
 has measure zero, and Theorem \ref{thm1} follows.

 \section{Proof of Proposition \ref{prop3}}\label{proof3}
 
Let $F$ be a smooth non-negative function with compact support in \eqref{inside} with $F({\bm \alpha}) \geq 1$ for $\bm \alpha \in \mathcal{D}^{\ast}$. We write $\| F \| := \int F(\alpha) d_{\text{hyp}}\alpha$. We continue to assume \eqref{NDelta}. It suffices to estimate 
 $$\int_{\mathcal{D}^{\ast}} F(\bm \alpha) \Big|\mathcal{G}_{\bm \alpha}(M, T) - \widehat{V^{\ast}}(0) \widehat{W}(0) \frac{M^2}{T} \Big|^2 d^{\ast}_{\text{{\rm hyp}}}{\bm \alpha}\ll_{\varepsilon} \frac{M^{4 - \frac{1}{44}+\varepsilon}}{T^2} + \frac{M^{2+\varepsilon}}{T} .$$
Opening the square, the proposition follows from the following two estimates
 \begin{equation}\label{first}
 \mathcal{I}_1(M, T) = \int_{\mathcal{D}^{\ast}}    F({\bm \alpha})\mathcal{G}_{\bm \alpha}(M, T) d^{\ast}_{\text{hyp}}{\bm \alpha}= \| F \|  \widehat{V^{\ast}}(0) \widehat{W}(0) \frac{M^2}{T} + O\Big( \frac{M^{3/2+\varepsilon}}{T} + 1\Big)
 \end{equation}
 and
   \begin{equation}\label{second}
   \begin{split}
  \mathcal{I}_2(M, T) &= \int_{\mathcal{D}^{\ast}}    F({\bm \alpha})|\mathcal{G}_{\bm \alpha}(M, T)|^2 d^{\ast}_{\text{hyp}}{\bm \alpha}\\
  &= \| F \|  \widehat{V^{\ast}}(0)^2 \widehat{W}(0)^2 \frac{M^4}{T^2} + O\Big(\frac{M^{4 - \frac{1}{44}+\varepsilon}}{T^2} + \frac{M^{2+\varepsilon}}{T} \Big).
  \end{split}
\end{equation} 
Before we proceed, we remark that without changing $\mathcal{G}_{\bm \alpha}(M, T)$ we can (and will) insert a redundant function $\psi(\textbf{x}/M,  \textbf{y}/M)$ into the sum, where $\psi$ a smooth function that is constantly 1 on some sufficiently large box in $\Bbb{R}^4$ (depending on the support of $F$).

\subsection{Proof of \eqref{first}}

 Fourier-inverting $V$ and $W$ we obtain
\begin{displaymath}
\begin{split}
 \mathcal{I}_1(M, T ) = &\frac{1}{4} \int_{\mathcal{D}^{\ast}} \int_{\Bbb{R}^2} F({\bm \alpha}) \sum_{(x_1, y_1) \not= \pm (x_2, y_2)}\psi\Big(\frac{\textbf{x}}{M},  \frac{\textbf{y}}{M}\Big)\widehat{V}(z) \widehat{W}(u) \\
 &e\Big(Tu\cdot \frac{q_{\bm \alpha}(x_1, y_1) - q_{\bm \alpha}(x_2, y_2)}{D({\bm \alpha}) }\Big)e\Big( z\frac{q_{\bm \alpha}(x_1, y_1)}{M^2 D({\bm \alpha})}\Big) du \, dz \frac{ d\alpha_1\, d\alpha_3}{\pi^3 D(\bm \alpha)^{3}}.
 \end{split}
 \end{displaymath}
 Changing variables $z \leftarrow z/D({\bm \alpha})$, $u \leftarrow  u/D({\bm \alpha})$ and integrating over ${\bm \alpha}$ we obtain
 \begin{displaymath}
\begin{split}
 \mathcal{I}_1(M, T ) = &\frac{1}{4 }   \int_{\Bbb{R}^2}   \sum_{(x_1, y_1) \not= \pm (x_2, y_2)}\psi\Big(\frac{\textbf{x}}{M},  \frac{\textbf{y}}{M}\Big)e\Big( T u ( x_1y_1 - x_2y_2)+z\frac{x_1y_1}{M^2}  \Big) \\
 & \widehat{G}\Big(T u (x_1^2 - x_2^2)  + z\frac{x_1^2}{M^2}, Tu(y_1^2 - y_2^2)+z \frac{y_1^2}{M^2} ; z, u\Big) du \, dz   \end{split}
 \end{displaymath}
where 
$$G(\alpha_1, \alpha_3;z, u) = \frac{F({\bm \alpha})}{\pi^3 D(\bm \alpha)}\widehat{V}(z D({\bm \alpha})) \widehat{W}(u D({\bm \alpha}))$$
and the Fourier transform is taken with respect to the first two variables. Note that $\widehat{G}$ is a Schwartz-class function in all variables. We change variables
\begin{equation}\label{change}
a_1 = y_1- y_2, \quad a_2 = y_1 + y_2, \quad b_1 = x_1 - x_2, \quad b_2 = x_1 + x_2,
\end{equation}
so that 
\begin{equation}\label{a}
a_1 \equiv a_2\, (\text{mod }2), \quad b_1 \equiv b_2\, (\text{mod } 2), \quad (a_1, b_1) \not= (0, 0) \not= (a_2, b_2).
\end{equation}
 We obtain
 \begin{equation}
 \label{I1}
  \mathcal{I}_1(M, T ) = \frac{1}{4} \left. \sum_{\textbf{a}, \textbf{b}}\right.^{\ast}  \tilde{\psi}\Big(\frac{\textbf{a}}{M},  \frac{\textbf{b}}{M}\Big) \mathcal{H}(\textbf{a}, \textbf{b})
  \end{equation}
 where the star indicates that the $\textbf{a}, \textbf{b}$-sum is subject to \eqref{a}, 
 $$ \tilde{\psi} (\textbf{a} ,   \textbf{b} ) = \psi\Big(\frac{b_1+b_2}{2}, \frac{b_2 - b_1}{2}, \frac{a_1+a_2}{2}, \frac{a_2-a_1}{2} \Big)$$
  and 
  \begin{equation}\label{defH}
\begin{split}
 \mathcal{H}(\textbf{a}, \textbf{b}) = &   \int_{\Bbb{R}^2}   e\Big(\frac{1}{2}Tu (a_1b_2 + a_2b_1) +\frac{z}{4M^2} (a_1+a_2)(b_1+b_2) \Big) \\
 & \widehat{G}\Big( Tu b_1b_2 + \frac{z}{4M^2}(b_1+b_2)^2, Tua_1a_2+ \frac{z}{4M^2} (a_1+a_2)^2; z, u\Big) du \, dz .  \end{split}
 \end{equation}
  Let 
 \begin{displaymath}
 \begin{split}
  &  C_{\max} := \max(|a_1|, |a_2|, |b_1|,  |b_2|), \quad C_{\min} := \min(|a_1|, |a_2|, |b_1|,  |b_2|),\\
   & P = \max(|a_1a_2|, |b_1b_2|, |a_1b_2 + a_2b_1|).
   \end{split}
   \end{displaymath}
 By \eqref{a} we have $C_{\max}, P \not= 0$.    Integrating by parts with respect to $u$ we find
\begin{equation}\label{H}
 \mathcal{H}(\textbf{a}, \textbf{b}) \ll_A     \frac{1}{T(|a_1a_2| + |b_1b_2|) + 1} \Big(1 + \frac{|a_1b_2 + a_2b_1|}{|a_1a_2| + |b_1b_2| + 1/T}\Big)^{-A} \ll \frac{1}{TP+1}
 \end{equation}
for any $A > 0$.   (If $\textbf{a}, \textbf{b}$ are integral vectors satisfying \eqref{a}, we have $TP + 1 \asymp TP$, but for arbitrary arguments it is important to keep the extra $+1$).  This argument shows that we can restrict the $u$-integral in \eqref{defH} to 
 \begin{equation}\label{u}
 T u P \ll M^{\varepsilon}
 \end{equation}
up to a negligible error.

 For some $1 < C < M$ to be determined later we first estimate the contribution of tuples $(\textbf{a}, \textbf{b})$ with $C_{\min} \leq C$ to \eqref{I1}. Assume without loss of generality that $C_{\min} = |a_1|$. If  $a_1= 0$, we distinguish the cases $b_2 \not= 0$ (in which case $b_1b_2 \not=0$ by \eqref{a}) and $b_2 = 0$ (in which case $\mathcal{H}(\textbf{a}, \textbf{b})$ is negligible by \eqref{H}). It is then easy to see that we get a contribution of 
$$\ll \frac{M(\log M)^2}{T}$$
to \eqref{I1}. 
If $|a_1| > 0$, then by \eqref{H} we obtain a contribution
 $$\sum_{0 \not= |a_1| \leq C}\sum_{ 0 \not=  a_2, b_1, b_2 \ll M } \frac{1}{T |b_1b_2|} \ll
\frac{M^{1+\varepsilon} C}{T}.$$
We now attach a smooth cut-off function  
$\Psi ( \textbf{a}/C,  \textbf{b}/C )$ 
to  the 
 $(\textbf{a}, \textbf{b})$-sum in \eqref{I1}
 where $\Psi$ is supported on $(-\infty, -1/2] \cup [1/2, \infty)$ in each variable and constantly 1 on $(-\infty, -1] \cup [1, \infty)$ in each variable. By the above remarks we thus obtain
 $$ \mathcal{I}_1(M, T ) = \frac{1}{4} \left. \sum_{\textbf{a}, \textbf{b}}\right.^{\ast} \Psi\Big(\frac{\textbf{a}}{M},  \frac{\textbf{b}}{M}\Big)  \tilde{\psi}\Big(\frac{\textbf{a}}{M},  \frac{\textbf{b}}{M}\Big) \mathcal{H}(\textbf{a}, \textbf{b})+ O\Big(\frac{M^{1+\varepsilon} C}{T}\Big).$$
 Recalling \eqref{u} and \eqref{H}, we  conclude
 \begin{equation*}
 \begin{split}
 \Big\|  \nabla\Big(\Psi\Big(\frac{\textbf{a}}{C}, \frac{\textbf{b}}{C} \Big)  \tilde{\psi}\Big(\frac{\textbf{a}}{M},  \frac{\textbf{b}}{M}\Big) \mathcal{H}(\textbf{a}, \textbf{b})\Big)\Big\| & \ll \int_{u \ll  M^{\varepsilon}/TP}  \Big(T|u|C_{\max} +\frac{C_{\max}}{M^2} \Big)   du  + \frac{1}{CTP} \\
 & \ll M^{\varepsilon}   \frac{C_{\max}}{P}  \frac{1}{TP} + \frac{1}{CTP} \ll \frac{M^{\varepsilon}}{CTP}. 
  \end{split}
 \end{equation*}
By  the Euler-MacLaurin formula  we obtain
 altogether
\begin{displaymath}
\begin{split}
\mathcal{I}_1(M, T) = \frac{1}{4}  \cdot \frac{1}{4} \int_{\Bbb{R}^4}  \Psi\Big(\frac{\textbf{a}}{C}, \frac{\textbf{b}}{C} \Big) \tilde{\psi}\Big(\frac{\textbf{a}}{M},  \frac{\textbf{b}}{M}\Big) \mathcal{H}(\textbf{a}, \textbf{b}) d(\textbf{a}, \textbf{b})+ O\Big( \frac{M^{2+\varepsilon}}{CT} + \frac{M^{1+\varepsilon}C}{T}\Big).
\end{split} 
\end{displaymath}
We can now remove  $\Psi(\textbf{a}/C,  \textbf{b}/C)$ in the same way as we introduced it at the cost of an error 
$$\ll \int_{|a_1| \leq C} \int_{|a_1| \ll a_2, b_1, b_2 \ll M} \frac{1}{T|b_1b_2|+1}d(\textbf{a}, \textbf{b}) \ll  \frac{M^{1+\varepsilon}C}{T}$$ 
by \eqref{H}. Choosing $C = M^{1/2}$, we obtain
\begin{equation}\label{approx}
\begin{split}
\mathcal{I}_1(M, T) = \frac{1}{4}  \cdot \frac{1}{4} \int_{\Bbb{R}^4}   \tilde{\psi}\Big(\frac{\textbf{a}}{M},  \frac{\textbf{b}}{M}\Big) \mathcal{H}(\textbf{a}, \textbf{b}) d(\textbf{a}, \textbf{b})+ O\Big( \frac{M^{3/2+\varepsilon}}{T}  \Big).
\end{split} 
\end{equation}
At this point we revert all transformations. We change variables $(\textbf{a}, \textbf{b})$ back to $(\textbf{x}, \textbf{y})$ using \eqref{change}, 
 undo the Fourier inversion with respect to ${\bm \alpha}$, the linear changes of variables of $(z, u)$, and the Fourier inversion with respect to  $u, z$. In this way the main term in \eqref{approx} equals 
  \begin{equation*} 
  \frac{1}{4}\int_{\mathcal{D}^{\ast}} \int_{\Bbb{R}^4}   F({\bm \alpha})  \psi\Big(\frac{\textbf{x}}{M},  \frac{\textbf{y}}{M}\Big) W\Big( T\cdot \frac{q_{\bm \alpha}(x_1, y_1) - q_{\bm \alpha}(x_2, y_2)}{D({\bm \alpha}) }\Big) V\Big( \frac{q_{\bm \alpha}(x_1, y_1)}{M^2 D({\bm \alpha})}\Big) d\textbf{x}\, d\textbf{y}  \,  d_{\text{hyp}}{\bm \alpha} .
  \end{equation*}
We drop $\psi$ because it is redundant. Let
$$x' = x +\frac{y}{2\alpha_1}, \quad x'' = \frac{x' \alpha_1^{1/2}}{(4\alpha_1\alpha_3 - 1)^{1/4}},  \quad y' = y \frac{(4\alpha_1\alpha_3 - 1)^{1/4}}{(4\alpha_1)^{1/2}}.$$  Then $$\frac{q_{\bm \alpha}(x, y)}{D({\bm \alpha})} = \pi \frac{\alpha_1 (x')^2 +  (4\alpha_1\alpha_3 - 1) y^2/(4\alpha_1 )}{(4\alpha_1\alpha_3 - 1)^{1/2} }= \pi \big((x'')^2 + (y')^2\big).$$ 
Changing variables, we obtain
 \begin{equation*} 
 \| F \|\int_{\Bbb{R}^4}     W\Big( T\pi  (x_1^2 + y_1^2 - x_2^2 - y^2_2) \Big) V\Big(\pi  \frac{(x_1^2 - y_1^2)}{M^2  }\Big) d\textbf{x}\, d\textbf{y}   .
  \end{equation*}
Changing to polar coordinates and recalling \eqref{vast}, we get
\begin{displaymath}
\begin{split}
 & \| F \| \int_{[0, \infty)^2}   W\Big( T  (r_1 - r_2) \Big) V\Big(  \frac{r_1}{M^2  }\Big)   dr_2\, dr_1\\
  &=   \|F \| \frac{M^2}{T}\int_{0}^{\infty} \int_{-\infty}^{M^2r_1/T}  W ( r_2 
) V (r_1 )   dr_2\, dr_1 =  \| F \| \frac{M^2}{T}\Big(\widehat{W}(0) \widehat{V^{\ast}}(0) + O\Big( \frac{T}{M^2}\Big)\Big). 
\end{split}
\end{displaymath}
Together with \eqref{approx} this establishes  \eqref{first}. 

\subsection{ Proof of \eqref{second}}
This follows to some extent the analysis in the proof of  \cite[Proposition 5]{ABR}. By Fourier inversion we have 
 \begin{equation*}
\begin{split}
\mathcal{I}_2(T, M)= &\frac{1}{16} \int_{\mathcal{D}^{\ast}}  F({\bm \alpha})\sum_{ \substack{(x_1, y_1) \not= \pm (x_2, y_2)\\ (x_3, y_3) \not= \pm (x_4, y_4)}}\psi\left( \frac{\textbf{x}}{M}, \frac{\textbf{y}}{M}\right) 
\int_{\Bbb{R}^4} \widehat{V}(z_1) \widehat{V}(z_2) \widehat{W}(u)\widehat{W}(v)\\
&e\left(Tu\frac{q_{{\bm \alpha}}(x_1, y_1) - q_{{\bm \alpha}}(x_2, y_2)}{D({\bm \alpha}) }\right) e\left(Tv\frac{q_{{\bm \alpha}}(x_3, y_3) - q_{{\bm \alpha}}(x_4, y_4)}{D({\bm \alpha}) }\right) \\
&e\left(\frac{z_1 q_{{\bm \alpha}}(x_1, y_1) - z_2 q_{{\bm \alpha}}(x_3, y_3)}{M^2 D({\bm \alpha}) }\right) \,  d(u, v, z_1, z_2)\,  \frac{d{\bm \alpha}}{\pi^3 D({\bm \alpha})^3},
 \end{split}
\end{equation*}
where $\psi$ is a similar redundant function as before. 
We change variables $ z_j\leftarrow z_j/D({\bm \alpha}) $, $u \leftarrow u/D({\bm \alpha}) $, $v \leftarrow v/D({\bm \alpha}) $, integrate over ${\bm \alpha}$ and introduce the variables
\begin{equation*}
\begin{split}
  & a_1 = y_1 - y_2, \quad a_2 = y_1 + y_2,\quad a_3 = x_3 - x_4,\quad a_4 =x_3 + x_4 ,\\
  &b_1 = x_1 - x_2, \quad b_2 = x_1 + x_2,\quad b_3 = y_3 - y_4,\quad b_4 = y_3 + y_4,\\
\end{split}
\end{equation*}
obtaining 
\begin{equation}\label{I2}
\mathcal{I}_2(T, M) = \frac{1}{16} \underset{\textbf{a}, \textbf{b} }{\left.\sum\right.^{\prime}} \tilde{\psi}\left( \frac{\textbf{a}}{M}, \frac{\textbf{b}}{M}\right)  \mathcal{H}(\textbf{a}, \textbf{b}) 
\end{equation}
 where $\sum^{\prime}$ indicates the conditions
\begin{equation}\label{congr}
a_1 \equiv a_2 \, (\text{mod } 2), \quad b_1 \equiv b_2 \, (\text{mod } 2),\quad a_3 \equiv a_4 \, (\text{mod } 2),\quad b_3 \equiv b_4 \, (\text{mod } 2)
 \end{equation}
 as well as
 \begin{equation}\label{ab}
  (a_1, b_1) \not= (0, 0) \not= (a_2, b_2), \quad  (a_3, b_3) \not= (0, 0) \not= (a_4, b_4).
  \end{equation}
 Moreover,  
 $$
 \tilde{\psi}(\textbf{a}, \textbf{b}) := \psi \Big ( \frac{b_1 + b_2}{2} , \frac{b_2 - b_1}{2}, \frac{a_3 + a_4}{2}, \frac{a_4 - a_3}{2}, \frac{a_1 + a_2}{2}, \frac{a_2 - a_1}{2}, \frac{b_3 + b_4}{2}, \frac{b_4 - b_3}{2} \Big ) = \psi(\textbf x, \textbf y)
 $$
and $ \mathcal{H}(\textbf{a}, \textbf{b}) $ is defined by 
 \begin{displaymath}
\begin{split}
    \int_{\Bbb{R}^4}  & e\Big(\frac{Tu}{2}(a_1b_2 + a_2 b_1) + \frac{Tv}{2}(a_3b_4 + a_4b_3) + \frac{z_1(a_1+a_2)(b_1+b_2)  - z_2 (a_3+a_4)(b_3+b_4) }{4M^2}  \Big) \\
 &  \widehat{G}\Big(Tub_1b_2 + Tva_3a_4 + \frac{z_1(b_1+b_2)^2  - z_2(a_3+ a_4)^2 }{4M^2} ,  \\
 & \quad  Tua_1a_2 + Tvb_3b_4 + \frac{z_1(a_1+a_2)^2 - z_2 (b_3 + b_4)^2}{4M^2}    ;z_1, z_2, u, v\Big)  \,  d(u, v, z_1, z_2)
 \end{split}
 \end{displaymath}
with 
 $$G(\alpha_1, \alpha_3; z_1, z_2, u, v) = \pi^{-3} D(\bm \alpha) F({\bm \alpha})  \widehat{V}(z_1 D({\bm \alpha}) ) \widehat{V}(z_2 D({\bm \alpha}) )  \widehat{W}(uD({\bm \alpha}) )\widehat{W}(vD({\bm \alpha}) ) $$
and the Fourier transform $\widehat{G}$ of $G$ is taken with respect to the first two variables $\alpha_1, \alpha_3$. We have 
\begin{equation}\label{schwartz}
\begin{split}
&\mathscr{D}\widehat{G}(U, V; z_1, z_2, u, v)\\
& \ll_{A, \mathscr{D}} \big((1+|U|)(1 + |V|)(1+|z_1|)(1+|z_2|)(1 + |u|)(1+|v|)\big)^{-A}
\end{split}
\end{equation}
for all $A > 0$ and any differential operator $\mathscr{D}$ with constant coefficients. 
Put 
\begin{displaymath}
\begin{split}
&P = \max(|a_1 a_2|, |a_3a_4|, |b_1b_2|, |b_3b_4|), \\
& \Delta = a_1a_2a_3a_4 - b_1b_2b_3b_4,\\
& \Delta_1 = a_1a_2b_3a_4 + a_1 a_2 a_3 b_4 - a_1 b_2b_3b_4-b_1a_2b_3b_4, \\
& \Delta_2 = a_1b_2a_3a_4 + b_1 a_2 a_3 a_4 - b_1 b_2a_3b_4-b_1b_2b_3a_4.
\end{split}
\end{displaymath}
We see immediately that the contribution of $P = 0$ to \eqref{I2} is negligible by \eqref{ab}, since in this case $a_1a_2 = b_1b_2 = a_3a_4 = b_3b_4 = 0$, but $a_1b_2 + a_2 b_1 \not= 0 \not=  a_3b_4 + a_4b_3$, so that  repeated partial integration in $u$ or $v$  saves as many powers of $T$ as we wish (and  $T \rightarrow \infty$  by \eqref{MT}). From now on we restrict to $P \not= 0$. If $\Delta\not = 0$, we change variables
$$u = \frac{a_3a_4 V - b_3b_4 U}{\Delta}, \quad v = \frac{a_1a_2 U - b_1b_2V}{\Delta}$$
to see that $\mathcal{H}(\textbf{a}, \textbf{b}) $ equals
\begin{equation}\label{H-var}
\begin{split}   \frac{1}{T^2| \Delta|}&\int_{\Bbb{R}^4}   e\Big(\frac{U\Delta_1}{2\Delta} + \frac{V\Delta_2}{2\Delta}   + \frac{z_1(a_1+a_2)(b_1+b_2) - z_2 (a_3+a_4)(b_3+b_4) }{4M^2}   \Big) \\
& \widehat{G}\Big(U + \frac{z_1(b_1+b_2)^2 - z_2(a_3+ a_4)^2 }{4M^2}   ,  V + \frac{z_1(a_1+a_2)^2 - z_2 (b_3 + b_4)^2 }{4M^2} ; \\
&\quad\quad z_1, z_2, \frac{a_3a_4 V - b_3b_4 U}{\Delta T}, \frac{a_1a_2 U - b_1b_2V}{\Delta T} \Big)   d(U, V, z_1, z_2). 
 \end{split}
 \end{equation}
By \eqref{schwartz} and repeated integration by parts we conclude
\begin{equation}\label{boundH}
\mathcal{H}(\textbf{a}, \textbf{b}) \ll_A \frac{1}{T^2(|\Delta| + P/T)} \Big(\Big(1 + \frac{|\Delta_1|}{|\Delta| + P/T}\Big)\Big(1 + \frac{|\Delta_2|}{|\Delta| + P/T}\Big)\Big)^{-A}
\end{equation}
for any $A \geq 0$ and all $\textbf{a}, \textbf{b} \ll M$ satisfying $P\not= 0$. This remains true for $\Delta = 0$, in which case we change variables
$$ua_1a_2 + v b_1b_2 = U$$
so that $ub_1b_2 + va_3a_4 = Ub_1b_2/a_1a_2$ and apply the same argument, cf.\  also \cite[(2.13)]{ABR}.  
We also conclude that
\begin{equation}\label{Hdiff}
\big\| \nabla\mathcal{H}(\textbf{a}, \textbf{b})\big \|  \ll \frac{1}{T^2|\Delta| }  \Big( \frac{1}{M}  +   \frac{M^3}{\Delta}\Big).
\end{equation}
We finally observe the trivial bound
\begin{equation}\label{Htriv}
  \mathcal{H}(\textbf{a}, \textbf{b}) \ll 1, 
  \end{equation}
valid for all real $\textbf{a}, \textbf{b}$ (even in the case $P=0$).  Fix $0 < \delta < 1/2$. We claim that the error from dropping terms in $\mathcal{I}_2(M, T)$ in \eqref{I2} with $\Delta \ll M^{4-\delta}$ is  small, more precisely (recall \eqref{boundH})
\begin{equation}\label{claim}
\begin{split}
& \frac{\#\{(\textbf{a}, \textbf{b})  \text{ satisfying } \eqref{ab} \mid \textbf{a}, \textbf{b} \ll M, \Delta \leq M^{4 - \delta}, \Delta_1, \Delta_2 \leq M^{\varepsilon}(|\Delta| + P/T)\} }{T^2(|\Delta| + P/T)} \\
&\ll M^{\varepsilon}\Big(\frac{M^{4-\delta/5}}{T^2} + \frac{M^2}{T}\Big). 
\end{split}
\end{equation}
 This is the analogue of \cite[Proposition 5]{ABR}. We postpone the proof  of \eqref{claim} to the next subsection. Let $\phi$ be a smooth function with support 
on $[1/2, \infty]$ that is 1 on $[1, \infty]$, and  write
\begin{equation*}
\Phi(\textbf{a}, \textbf{b}) :=  \tilde{\psi}\left( \frac{\textbf{a}}{M}, \frac{\textbf{b}}{M}\right)  \phi\left(\frac{|\Delta|}{M^{4-\delta}}\right) ,
\end{equation*} 
so that 
$$\mathcal{I}_2(T, M) = \frac{1}{16} \underset{\textbf{a}, \textbf{b} }{\left.\sum\right.^{\prime}}\Phi(\textbf{a}, \textbf{b}) \mathcal{H}(\textbf{a}, \textbf{b}) + O\Big(M^{\varepsilon}\Big(\frac{M^{4-\delta/5}}{T^2} + \frac{M^2}{T}\Big)\Big).$$ 
Note that the condition \eqref{ab} is now void. From \eqref{boundH} and \eqref{Hdiff} we conclude
$$\big\|\nabla \Phi(\textbf{a}, \textbf{b}) \mathcal{H}(\textbf{a}, \textbf{b}) \big \| \ll \frac{1}{T^2|\Delta| } \Big( \frac{1}{M^{1-\delta}}  +   \frac{M^3}{\Delta}   \Big)\ll \frac{1}{T^2 M^{5-2\delta}}.  $$
By the Euler-MacLaurin formula we conclude
$$ \underset{\textbf{a}, \textbf{b} }{\left.\sum\right.^{\prime}}\Phi(\textbf{a}, \textbf{b}) \mathcal{H}(\textbf{a}, \textbf{b}) = \frac{1}{16} \int_{\Bbb{R}^4} \int_{\Bbb{R}^4} \Phi(\textbf{a}, \textbf{b}) \mathcal{H}(\textbf{a}, \textbf{b})  d\textbf{a} \, d\textbf{b} + O\Big(  \frac{M^{3+2\delta}}{T^2}\Big),$$
where the factor $1/16$ comes from  the congruences \eqref{congr}. 
We now re-insert the contribution of the terms $\Delta \ll M^{4-\delta}$ into the integral by dropping the cut-off function $\phi(|\Delta|/M^{4-\delta})$, and we claim that this introduces an error of at most
\begin{equation}\label{claim2}
 \int_{\Delta \ll M^{4-\delta}} \tilde{\psi}\left( \frac{\textbf{a}}{M}, \frac{\textbf{b}}{M}\right)\mathcal{H}(\textbf{a}, \textbf{b}) d(\textbf{a}, \textbf{b})  \ll   \frac{M^{4-\delta/5+\varepsilon}}{T^2}, 
\end{equation}
which is already present. Again we postpone the proof and  revert all steps as in the previous subsection and in the end of \cite[Section 2]{ABR}, namely the change of variables $(\textbf{x}, \textbf{y}) \mapsto (\textbf{a}, \textbf{b})$, the integration over $\bm \alpha$ and the Fourier inversions. In this way we finally obtain that $\mathcal{I}_2(M, T)$ equals
\begin{displaymath}
\begin{split}
 \frac{1}{16} \int_{\mathcal{D}^{\ast}} &F({\bm \alpha}) \int_{\Bbb{R}^4} \int_{\Bbb{R}^4}   
\psi\left(\frac{\textbf{x}}{M}, \frac{\textbf{y}}{M}\right) W\left(\frac{q_{{\bm \alpha}}(x_1, y_1) - q_{{\bm \alpha}}(x_2, y_2)}{D({\bm \alpha}) }\right)\\
&W\left(\frac{q_{{\bm \alpha}}(x_3, y_3) - q_{{\bm \alpha}}(x_4, y_4)}{D({\bm \alpha}) }\right)V\left(\frac{q_{{\bm \alpha}}(x_1, y_1)}{M^2 D({\bm \alpha}) }\right) V\left(\frac{q_{{\bm \alpha}}(x_3, y_3)}{M^2 D({\bm \alpha}) }\right) \, d\textbf{x} \, d\textbf{y} \, d_{\text{hyp}}{\bm \alpha} \\
&+ O\Big(M^{\varepsilon}\Big(\frac{M^{4-\delta/5}}{T^2} + \frac{M^2}{T} +  \frac{M^{3+2\delta}}{T^2}\Big)\Big). 
\end{split}
\end{displaymath} 
Here we can drop the function $\psi(\textbf{x}/M, \textbf{y}/M)$ because it is redundant. We choose   $\delta = 5/11$. By the same change of variables as in the previous subsection we obtain
\begin{displaymath}
\begin{split}
\mathcal{I}_2(T, M) & =  \| F \| \int_{\Bbb{R}^8} W\Big( T\pi  (x_1^2 + y_1^2 - x_2^2 - y^2_2) \Big)W\Big( T\pi  (x_3^2 + y_3^2 - x_4^2 - y^4_2) \Big) \\
&\quad\quad V\Big(\pi  \frac{(x_1^2 - y_1^2)}{M^2  }\Big) V\Big(\pi  \frac{(x_3^2 - y_3^2)}{M^2  }\Big) d\textbf{x}\, d\textbf{y}   + O\Big( \frac{M^{43/11+\varepsilon}}{T^2} + \frac{M^{2+\varepsilon}}{T}  \Big).
\end{split}
\end{displaymath}
The main term equals
\begin{displaymath}
\begin{split}
& \| F \| \frac{M^4}{T^2}\int_{0}^{\infty} \int_{-\infty}^{M^2r_1/T} \int_{0}^{\infty} \int_{-\infty}^{M^2r_3/T}  W ( r_2 
) V (r_1 ) W(r_4) V(r_3) dr_4\, dr_3 \,  dr_2\, dr_1\\
=&\| F \| \frac{M^4}{T^2} \Big(\widehat{W}(0) \widehat{V^{\ast}}(0) + O\Big(\frac{T}{M^2}\Big)\Big)^2,
\end{split}
\end{displaymath}
and \eqref{second} follows.  It remains to prove \eqref{claim} and \eqref{claim2} to which the following two subsections are devoted.  

\subsection{Proof of \eqref{claim}}

We use the notation $X \preccurlyeq Y$ to mean $X \ll M^{\varepsilon} Y$. 
We put all variables into dyadic intervals and suppose that 
 $$A _1 \leq |a_1| \leq  2A_1, \ldots, A_4 \leq |a_4| \leq 2 A_4, B_1 \leq |b_1| \leq  2B_1, \ldots, B_4 \leq |b_4| \leq 2 B_4$$ with $0 \leq A_1, \ldots, B_4 \ll M$. We also assume 
 \begin{equation}\label{D}
    D \leq |\Delta| \leq 2D \ll M^{4-\delta}.
    \end{equation}
We now   count the number $\mathcal{N}(\textbf{A}, \textbf{B}, D)$ of 8-tuples $(\textbf{a}, \textbf{b})$ subject to these size conditions as well as \eqref{ab} and
\begin{equation}\label{delta}
\Delta_1, \Delta_2 \preccurlyeq |\Delta| + P/T .
\end{equation}
We start with some degenerate cases and denote by $\mathcal{N}_0(\textbf{A}, \textbf{B}, D)$ the contribution of $A_1\cdots A_4B_1\cdots B_4D = 0$. 

Let us first assume that some of the $\textbf{a}$-variables vanish, but none of the $\textbf{b}$-variables. Then $D \not = 0$, and by a divisor bound, this contribution is $\preccurlyeq D M^3$. If some $\textbf{a}$-variable vanishes, say $a_1$, and also some $\textbf{b}$-variable vanishes (which cannot be $b_1$), say $b_2$, then $|b_1a_2| \geq 1$, $\Delta = 0$, $P = \max(|a_3a_4|, |b_3b_4|)$, $\Delta_1 = b_1a_2b_3b_4$, $\Delta_2 = b_1a_2a_3a_4$, which is impossible by \eqref{MT} and \eqref{delta}. 

So from now on we will assume $a_1\cdots a_4b_1\cdots b_4 \not= 0$. Let us next assume $D= 0$, i.e.\ $a_1\cdots a_4 = b_1\cdots b_4$. By a divisor bound  contributes $\preccurlyeq PM^2$. So from now on we assume $D\not= 0$, 
 and we have shown \begin{equation}\label{firsterror}
\mathcal{N}_0(\textbf{A}, \textbf{B}, D) \preccurlyeq D M^3 + PM^2. 
\end{equation}
Let $\mathcal{N}_{\ast}(\textbf{A}, \textbf{B}, D)$ denote the contribution with $A_1\cdots A_4B_1\cdots B_4D \not= 0$ and let us write  $ \min(A_1, \ldots, A_4, B_1, \ldots B_4) = M^{1-\eta}$, say, with $0 < \eta < 1$.   We have trivially
\begin{equation*}
\mathcal{N}_{\ast}(\textbf{A}, \textbf{B}, D) \preccurlyeq D \min(A_1A_2A_3A_4, B_1B_2B_3B_4) \ll DM^{4-\eta} .
\end{equation*}

We consider another degenerate case, namely $a_2a_4 = b_2b_4$, In this case let $d = a_1a_3 - b_1b_3 = \Delta/(b_2b_4) \not=0$. We have $\preccurlyeq DM^2$ choices for $(d, b_2, b_4, a_1, a_3)$, and then $b_1, b_3$ are determined by a divisor argument. This can be absorbed in the existing count \eqref{firsterror}, so that from now we assume $a_2a_4\not=b_2b_4$ and similarly $a_2a_3 \not= b_2b_3$.

Substituting $  b_1 =(a_1a_2a_3a_4 - \Delta)/(b_2b_3b_4)$ into the definition of $\Delta_1$, we obtain
$$\Delta_1 = \frac{-a_1}{b_2}(a_2a_3 -b_2b_3)(a_2a_4 - b_2b_4) + \frac{a_2}{b_2} \Delta  $$
so
\begin{equation}\label{subst}
(a_2a_3 -b_2b_3)(a_2a_4 - b_2b_4) \preccurlyeq \frac{B_2}{A_1}\Big(D + \frac{P}{T}\Big) + \frac{A_2}{A_1} \ll \Big(D + \frac{P}{T}\Big)M^{\eta}.
\end{equation}
By a standard lattice point argument \cite[p.\ 200-201]{Sa} the number of such (non-zero) 6-tuples is
$$\preccurlyeq (D + P/T)M^{2+\eta}.$$
Having fixed $a_2, a_3, a_4, b_2, b_3, b_4$, we are left with pairs $(a_1, b_1)$ satisfying
$$a_1 = \frac{b_1b_2b_3b_4}{a_2a_3a_4} + O\Big(\frac{D}{A_2A_3A_4}\Big)$$
so that we obtain in total
\begin{displaymath}
\begin{split}
\mathcal{N}(\textbf{A}, \textbf{B}, D)& \preccurlyeq D M^3 + PM^2 + \min\Big(DM^{4-\eta},\Big(D + \frac{P}{T}\Big)M^{3+\eta} \Big(\frac{D}{M^{3-3\eta}} + 1\Big)\Big)\\
&\leq D M^3 + PM^2 +\Big(D + \frac{P}{T}\Big)\big(M^{7/2} +  D^{1/5} M^{16/5}\big)\\
& \ll \Big(D + \frac{P}{T}\Big) \big(M^ {4 - \delta/5} + M^2T\big)
\end{split}
\end{displaymath}
by \eqref{D} (and since $\delta < 1/2$). The quantity on the left hand side of  \eqref{claim} is $$\preccurlyeq \max_{\textbf{A}, \textbf{B}, D} \frac{\mathcal{N}(\textbf{A}, \textbf{B}, D)}{T^2 (D + P/T)}$$
and so \eqref{claim} follows.

 \subsection{ Proof of \eqref{claim2}}
 Here we must estimate
 $$  \int_{\Delta \ll M^{4-\delta}} \tilde{\psi}\left( \frac{\textbf{a}}{M}, \frac{\textbf{b}}{M}\right)\mathcal{H}(\textbf{a}, \textbf{b}) d(\textbf{a}, \textbf{b}) $$
  based on the bounds \eqref{boundH} and \eqref{Htriv}. This is the (much simpler) continuous analogue of the previous subsection. We put all variables into dyadic intervals $A_j \leq |a_j| \leq 2A_j$, $B_j \leq |b_j| \leq 2B_j$, $D \leq |\Delta| \leq 2D$, $\Delta_1, \Delta_2 \preccurlyeq \Delta  + P$. The numbers $A_j, B_j, D$ run through logarithmically many positive and negative powers of 2 and are bounded by $M^{-100} \leq A_j, B_j  \ll M$, $M^{-100} \leq D \ll M^{4-\delta}$.  (If one of these is $\ll M^{-100}$, the coarsest trivial estimates suffice - this is the only point where \eqref{Htriv} is needed.)  We call the corresponding set $\mathcal{S}(\textbf{A}, \textbf{B}, D)$ and as before we write $ \min(A_1, \ldots, A_4, B_1, \ldots B_4) = M^{1-\eta}$. We have the trivial bound
$$\text{vol}(\mathcal{S}(\textbf{A}, \textbf{B}, D)) \ll DM^{4-\eta+\varepsilon}$$
  On the other hand, as in \eqref{subst} we see that the integration condition implies
  $$(a_2a_3 -b_2b_3)(a_2a_4 - b_2b_4) \preccurlyeq \frac{B_2}{A_1}\Big(D + \frac{P}{T}\Big) + \frac{A_2}{A_1} \ll \Big(D + \frac{P}{T}\Big)M^{\eta}$$
An easy computation shows that the volume of such 6-tuples is $\ll (D + P/T)M^{2+\eta + \varepsilon}$, so that we obtain the alternative bound
$$\text{vol}(\mathcal{S}(\textbf{A}, \textbf{B}, D)) \ll \Big(D + \frac{P}{T}\Big)M^{2+\eta + \varepsilon} \cdot M \cdot \frac{D}{M^{3-3\eta}}\ll \Big(D + \frac{P}{T}\Big)M^{4+ 4\eta - \delta+\varepsilon}.$$
Combining the two bounds we obtain   
 $$\text{vol}(\mathcal{S}(\textbf{A}, \textbf{B}, D)) \ll (D + P/T)M^{4 - \delta/5+\varepsilon}$$
 and \eqref{claim2} follows from this and \eqref{boundH}. 
 
  \section{Proof of Theorem \ref{thm2}}
  
 As in the proof of Theorem \ref{thm1} we can specialize $\alpha_3 = 1$, and we can restrict ourselves to a compact set of $\alpha \in \mathcal{R}_0 \subseteq \Bbb{R}_{> 0}$.  Thus our quadratic forms have the shape $q_{\alpha}(m, n) = \alpha m^2 + n^2$, $m > 0$, $n \geq 0$  with $\alpha \asymp 1$. 
By the same argument as in the proof of Theorem \ref{thm1} it suffices to show that the measure of $\alpha \in \mathcal{R}_0$ such that
\begin{equation}\label{suchthat}
 \big| \#\{\Lambda_i, \Lambda_j \leq N \mid 0 \leq \Lambda_j - \Lambda_i \leq \Delta, i \not= j \} - N\Delta \big| > \delta N\Delta
 \end{equation}
is $ \ll_{\delta, \eta} 1/(\log N)^{2.5}$ (the exponent has to be larger than 2 for the Borel-Cantelli argument to work), uniformly in the range \eqref{NDelta1}.  
We define 
\begin{equation}\label{var}
t_1 = m_1-m_2, \quad t_2 = m_1+m_2, \quad t_3 = n_2 - n_1, \quad t_4 = n_2 + n_1,
\end{equation}
so that $\#\{\Lambda_i, \Lambda_j \leq N \mid 0 \leq \Lambda_j - \Lambda_i \leq \Delta, i \not= j \}$ equals the cardinality of 
 all quadruples $(t_1, t_2, t_3, t_4) \in \mathcal{T}_{\alpha}(N, \Delta)$ where $\mathcal{T}_{\alpha}(N, \Delta)$ is defined by  
 \begin{displaymath}
\begin{split}
  & (t_1, t_3) \not= (0, 0), \quad   t_1 \equiv t_2\, (\text{mod } 2), \quad t_3 \equiv t_4\, (\text{mod }2),  \quad  t_2  > |t_1|, \quad t_4  \geq |t_3|, \\
  & \alpha \Big( \frac{t_1\pm t_2}{2}\Big) ^2 + \Big( \frac{t_3\pm t_4}{2}\Big)^2 \leq \frac{4\sqrt{\alpha}N}{\pi}, \quad 0 \leq \alpha t_1t_2 - t_3 t_4 \leq \frac{4\sqrt{\alpha} \Delta}{\pi}.
\end{split}
\end{displaymath}
For $\rho > 0$ let  
$$S_{\rho}(N) = \Big\{n \in \Bbb{N} \mid n \text{ has a prime divisor in } \big(\exp((\log N)^{\rho}), \exp((\log N)^{1-\rho})\big]\Big\}.$$
Let $\mathcal{T}^{\rho, A}_{\alpha}(N, \Delta)$ be the set of $(t_1, \ldots, t_4) \in \mathcal{T}_{\alpha}(N, \Delta)$ such that 
\begin{itemize}
\item at least one of $t_1, \ldots, t_3$ is in $S_{\rho}(N)$ (we make no assumption on $t_4$) and 
\item all $t_j$ satisfy $|t_j| \geq N^{1/2} (\log N)^{-A}$. 
\end{itemize}
We claim that the contribution of $(t_1, \ldots, t_4) \in \mathcal{T}_{\alpha}(N, \Delta) \setminus 
\mathcal{T}^{\rho, A}_{\alpha}(N, \Delta)$ is negligible. 

\begin{prop}\label{prop4} Let $A > 10$ and $\rho < 1/20$. Then the measure of all $\alpha \in \mathcal{R}_0$ such that
$$\# \big(\mathcal{T}_{\alpha}(N, \Delta) \setminus 
\mathcal{T}^{\rho, A}_{\alpha}(N, \Delta)\big)  \geq \frac{1}{2}\Delta N$$
is $\ll (\log N)^{-2.5}$. 
\end{prop}

We postpone the proof to Section \ref{proof4}.  

By symmetry  we can  assume that $t_1 \geq 0$ so that automatically $t_3 \geq 0$. Moreover, we drop the conditions $t_2 >  t_1$, $t_4 \geq t_3$.  In this way we see that $\#  
\mathcal{T}^{\rho, A}_{\alpha}(N, \Delta)$ is 1/2 times the cardinality of all $(t_1, \ldots, t_4)$ such that
\begin{displaymath}
\begin{split}
  &  t_1 \equiv t_2\, (\text{mod } 2), \quad t_3 \equiv t_4\, (\text{mod }2), \quad t_j \geq N^{1/2} (\log N)^{-A}\\
  & \alpha \Big( \frac{t_1\pm t_2}{2}\Big) ^2 + \Big( \frac{t_3\pm t_4}{2}\Big)^2 \leq \frac{4\sqrt{\alpha}N}{\pi}, \quad 0 \leq \alpha t_1t_2 - t_3 t_4 \leq \frac{4\sqrt{\alpha} \Delta}{\pi}
\end{split}
\end{displaymath}
and at least one of the $t_1, t_2, t_3$ is in $S_{\rho}(N)$. 

It is convenient to count a slightly smaller and a slightly large quantity. Fix some $\delta' > 0$. We restrict the $t_i$ to boxes $ D_i < t_i \leq (1 + \delta' )  D_i$, where the $D_i$ run through $(\log\log X)^{O_{\delta'}(1)}$ values of a sequence of the type $D(1+\delta')^j$. We can and will assume 
\begin{equation}\label{d1d2d3d4}
\frac{1}{2} \leq \frac{\alpha D_1D_2}{D_3D_4} \leq 2.
\end{equation}

For the smaller count   we only consider tuples $(D_1, \ldots, D_4)$ satisfying \eqref{d1d2d3d4} such that the boxes lie completely inside the region 
\begin{equation}\label{region}
\alpha \Big( \frac{t_1\pm t_2}{2}\Big) ^2 + \Big( \frac{t_3\pm t_4}{2}\Big)^2 \leq \frac{4\sqrt{\alpha}X}{\pi}, \quad t_1, \ldots, t_4 \geq 0.
\end{equation}
We call the collection of such quadruples $\mathcal{D}_-$. 
For $0 \leq \theta \leq 1$ we note that 
$$
\log\Big( \frac{t_3t_4}{t_1t_2} + \theta \frac{4\sqrt{\alpha}\Delta}{t_1t_2\pi}\Big) = \log \frac{t_3t_4}{t_1t_2} + \theta \frac{4\sqrt{\alpha}\Delta}{t_3t_4\pi} + O\Big( \frac{\theta^2\Delta^2}{D_3^2D_4^2}\Big) .$$
Therefore we sharpen the inequality 
\begin{equation}\label{inequ}
  0 \leq \alpha t_1t_2 - t_3 t_4 \leq \frac{4\sqrt{\alpha} \Delta}{\pi}
 \end{equation} 
   to 
$$0 \leq \log \alpha - \log  \frac{t_3t_4}{t_1t_2} \leq (1 - 3\delta')\frac{4 \sqrt{\alpha}\Delta}{\pi D_3D_4},$$
and we detect this with a smooth weight function 
$$W_{-} \Big( \Big(\log \alpha - \log  \frac{t_3t_4}{t_1t_2}\Big) \frac{\pi D_3D_4}{4\sqrt{\alpha}\Delta}\Big).$$
 where $W_-$ is constantly 1 on $[\delta', 1 - 4\delta']$ and vanishes outside $[0, 1-3\delta']$.

For the larger count we consider boxes satisfying \eqref{d1d2d3d4}  that intersect the part of \eqref{region} where $t_i \geq N^{1/2} (\log N)^{-A}$, calling this larger collection of quadruples $\mathcal{D}_+$. We relax \eqref{inequ} to  
$$0 \leq \log \alpha - \log  \frac{t_3t_4}{t_1t_2} \leq (1 + \delta')\frac{4 \sqrt{\alpha}\Delta}{\pi D_3D_4}$$
and detect it with a factor $$W_{+} \Big( \Big(\log \alpha - \log  \frac{t_3t_4}{t_1t_2}\Big) \frac{\pi D_3D_4}{4\sqrt{\alpha}\Delta}\Big).$$
 where $W_+$  is constantly 1 on $[0, 1 ]$ and vanishes outside $[-\delta', 1+\delta']$. Eventually we will choose $\delta' = c\delta$ for some sufficiently small  $c$. For notational simplicity we write
 $$D_j' = D_j(1 + \delta').$$
 We summarize that
\begin{equation*}
\begin{split}
\#  
\mathcal{T}^{\rho, A}_{\alpha}(N, \Delta)&\leq  \frac{1}{2}\sum_{(D_1, \ldots, D_4) \in \mathcal{D}_{+}} \sum_{\substack{t_i \in (D_i,  D_i']\\\text{some } t_1, t_2, t_3\in S_{\rho}(N)\\ 2 \mid t_1 - t_2, t_3 - t_4}} W_{+} \Big( \Big(\log \alpha - \log  \frac{t_3t_4}{t_1t_2}\Big) \frac{\pi D_3D_4}{4\sqrt{\alpha} \Delta}\Big)\\
& = \frac{1}{2} \sum_{(D_1, \ldots, D_4) \in \mathcal{D}_{+}} \frac{4\sqrt{\alpha}\Delta}{\pi D_3D_4}\int_{\Bbb{R}} \widehat{W}_{+}\Big(\frac{4\sqrt{\alpha}\Delta y}{\pi D_3D_4}\Big) \sum_{\substack{t_i \in (D_i, D'_i]\\ \text{some } t_1, t_2, t_3  \in S_{\rho}(N)\\ 2 \mid t_1 - t_2, t_3 - t_4}} \Big( \frac{\alpha t_1t_2}{t_3t_4}\Big)^{2 \pi i y} dy,
\end{split}
\end{equation*}
and we have a similar lower bound where the subscripts $+$ are replaced with the subscripts $-$. We consider only the $+$ case, the other one being identical. 
We extract the main term from small values of $y$ in the integral. Let $V$ be a smooth non-negative function that is 1 on $[-1, 1]$ and vanishes for $|x| > 2$. Let $B > 3A$,  define
\begin{equation}\label{defY}
Y :=   \frac{ Y' D_3D_4}{N}, \quad Y' = (\log N)^B
\end{equation}
and decompose the previous  $y$-integral as $I_1(\alpha) + I_2(\alpha)$ where $I_1(\alpha)$ contains the factor $V(y/Y)$ and $I_2(\alpha)$ contains the factor $1 - V(y/Y)$. 
We claim
\begin{prop}\label{prop5} We have 
$$\frac{1}{2} \sum_{(D_1, \ldots, D_4) \in \mathcal{D}_{+}} \frac{4\sqrt{\alpha}\Delta}{\pi D_3D_4} I_1(\alpha)= (1 + O(\delta')) \Delta N$$
where the implied constant   depends only on $\rho, A, B$ and $\mathcal{R}_0$. 
\end{prop}

\begin{prop}\label{prop6} For suitable choices of $A, B$ we have 
$$ \int_{\mathcal{R}_0} \Big|\frac{1}{2} \sum_{(D_1, \ldots, D_4) \in \mathcal{D}_{+}} \frac{4\sqrt{\alpha}\Delta}{\pi D_3D_4} I_2(\alpha)  \Big|^2 \frac{d\alpha}{\alpha} \ll \Delta^2 N^2 (\log N)^{-50} $$
where the implied constant  depends only on $\rho $ and $\mathcal{R}_0$. 
\end{prop}

We postpone the proofs to Sections \ref{proof5} and \ref{proof6} respectively. From these two propositions we obtain that the measure of $\alpha \in \mathcal{R}_0$ such that $$\big|\#  
\mathcal{T}^{\rho, A}_{\alpha}(N, \Delta) - \Delta N\big| \gg \delta' \Delta N$$
is $O((\log N)^{-50})$. We choose $\delta' = \delta/c$ where $c$ is sufficiently large in terms of the implied constant in Propositon \ref{prop5}. Together with Proposition \ref{prop4} we see that the measure of $\alpha \in \mathcal{R}_0$ satisfying \eqref{suchthat} is indeed $\ll (\log N)^{-2.5}$ which completes the proof.

 \section{Proof of Proposition \ref{prop4}}\label{proof4}
 
 We recall that $\mathcal{T}_{\alpha}(N, \Delta) \setminus 
\mathcal{T}^{\rho, A}_{\alpha}(N, \Delta)$ is contained in the  set of all $(t_1, \ldots, t_4)$ such that 
  \begin{equation}\label{general}
\begin{split}
 & (t_1, t_3) \not= (0, 0), \quad t_2 > |t_1|, \quad t_4 \geq |t_3|, \\
  &  \alpha \Big( \frac{t_1\pm t_2}{2}\Big) ^2 + \Big( \frac{t_3\pm t_4}{2}\Big)^2 \ll N, \quad 0 \leq \alpha t_1t_2 - t_3 t_4 \ll \Delta
\end{split}
\end{equation}
and none of $t_1, t_2, t_3 $ is in $S_{\delta}(N)$ or at least one $t_j$ satisfies $|t_j|  \leq N^{1/2} (\log N)^{-A}$. 

Let us first consider the contribution $\mathcal{T}_1(\alpha)$ of those quadruples, where at least one $t_i$, say $t_1$, satisfies $|t_1| \leq N^{1/2} (\log N)^{-A}$. We observe directly that the contribution of $t_1t_2t_3t_4 = 0$  is  bounded by $O(N^{1/2})$ if $\Delta \gg 1$ and vanishes otherwise, so it is $\ll N^{1/2} \Delta$. We  assume from now on that all $t_j$ are non-zero, and without loss of generality positive. 
 We localize each $t_i$ in a dyadic interval of the shape  $(D_i, 4D_i]$; there are at most $O((\log N)^4)$ such boxes, and we must have 
 \begin{equation}\label{DD}
    D_1D_2 \asymp D_3 D_4 \ll N (\log N)^{-A}
    \end{equation}
     and moreover
$$\Big|\log \alpha + \log \frac{t_1t_2}{t_3t_4} \Big|  \ll \frac{\Delta}{D_3D_4}.$$
If $W$ denotes a suitable non-negative smooth compactly supported function, we conclude that 
$$\mathcal{T}_1(\alpha) \ll  \sum_{D_1, D_2, D_3, D_4}\sum_{t_i \in (D_i, 4D_i]}W\Big(\frac{D_3D_4}{\Delta} \Big(\log \alpha + \log \frac{t_1t_2}{t_3t_4} \Big)\Big) + N^{1/2} \Delta$$
where the outer sum runs over powers of 2 subject to \eqref{DD}. Therefore the measure of $\alpha \in \mathcal{R}_0$ such that $\mathcal{T}_1(\alpha) > \Delta N/2$ is at most 
$$ \int_{\mathcal{R}_0} \frac{1}{\Delta N} \sum_{D_1, D_2, D_3, D_4}\sum_{t_i \in (D_i, 4D_i]}W\Big(\frac{D_3D_4}{\Delta} \Big(\log \alpha + \log \frac{t_1t_2}{t_3t_4} \Big)\Big) \frac{d\alpha}{\alpha}.$$
We write $\beta = \log \alpha$, insert the weight $e^{-\beta^2/2}$ and extend the integration to all of $\Bbb{R}$. By Fourier inversion we bound the previous display by
$$ \sum_{D_1, D_2, D_3, D_4} \frac{1}{\Delta N} \frac{\Delta}{D_3D_4} \int_{\Bbb{R}} e^{-u^2/2} \widehat{W}\Big( \frac{\Delta u}{D_3D_4}\Big) F(u, D_1)F(u, D_2)F(-u, D_3)F(-u, D_4)du$$
where
$$F(u, D) = \sum_{D < t \leq 4D} \frac{1}{t^{2\pi i u}}.$$ 
By trivial estimates and \eqref{DD}, the above is 
$$\ll   \sum_{D_1, D_2, D_3, D_4}  \frac{D_1D_2}{N} \ll (\log N)^{4-A}.$$

Let now $\mathcal{T}_2(\alpha)$ denote the set of all quadruples satisfying \eqref{general} with $|t_j| \geq N^{1/2} (\log N)^{-A}$, but none of  $t_1, t_2, t_3 $ is in $S_{\rho}(N)$. Again we localize each $t_i$ in a dyadic interval of the shape  $(D_i, 4D_i]$; by our current assumptions, there are at most $O_A((\log\log N)^4)$ such boxes. By the same argument as before, the measure of $\alpha \in \mathcal{R}_0$ such that $\mathcal{T}_2(\alpha) > \Delta N/2$ is at most 
\begin{equation}\label{T2}
\begin{split}
 \sum_{D_1, D_2, D_3, D_4} \frac{1}{\Delta N} \frac{\Delta}{D_3D_4} &\int_{\Bbb{R}} e^{-u^2/2} \widehat{W}\Big( \frac{\Delta u}{D_3D_4}\Big) \\
 &F_{\rho, N}(u, D_1)F_{\rho, N}(u, D_2)F_{\rho, N}(-u, D_3)F(-u, D_4)du
 \end{split}
 \end{equation}
where 
 $$F_{\rho, N}(u, D) = \sum_{\substack{D < t \leq 4D\\ t \not\in S_{\rho}(N)}} \frac{1}{t^{2\pi i u}}.$$ 
At this point we invoke \cite[Corollary 1.2(iii)]{We} with $y = \exp((\log N)^{1-\rho})$, $z = \exp((\log N)^{\rho})$, $u = (\log N)^{\rho}$, $r = (\log N)^{2\rho - 1}$ in the form
\begin{equation}\label{wein}
  \sum_{\substack{n \leq D\\ n \not\in S_{\rho}(N)}} 1  = D P_N+  O\big( D \exp(- (\log N)^{\rho/2})\big)
  \end{equation}
where
$$P_N = \prod_{\exp((\log N)^{\rho}) < p \leq \exp((\log N)^{1-\rho})} \Big(1 - \frac{1}{p}\Big) \ll (\log N)^{2\rho - 1},$$
 uniformly in $N^{1/3} \ll D \ll N^{1/2}$, say. By trivial bounds, \eqref{T2} is at most
 $$\ll  (\log \log N)^4 (\log N)^{ 6\rho-3},$$
 and the proof is complete. 
 
 \section{Proof of Proposition \ref{prop5}}\label{proof5}
 
 The aim of this section is an asymptotic evaluation of 
\begin{equation}\label{eval}
 \frac{1}{2} \sum_{(D_1, \ldots, D_4) \in \mathcal{D}_{+}} \frac{4\sqrt{\alpha}\Delta}{\pi D_3D_4}\int_{\Bbb{R}} \widehat{W}_{+}\Big(\frac{4\sqrt{\alpha}\Delta y}{\pi D_3D_4}\Big)V\Big( \frac{y}{Y}\Big) \sum_{\substack{t_i \in (D_i, D'_i]\\ \text{some } t_i  \in S_{\rho}(N)\\ 2 \mid t_1 - t_2, t_3 - t_4}} \Big( \frac{\alpha t_1t_2}{t_3t_4}\Big)^{2 \pi i y} dy
 \end{equation}
with $Y$ as in \eqref{defY}. By a Taylor argument we have
 $$\widehat{W}_+\Big(\frac{4\sqrt{\alpha}\Delta y}{\pi D_3D_4}\Big) =  \widehat{W}_+(0) + O_{\delta'}\Big( \frac{\Delta Y}{D_3D_4}\Big).$$
(The definition of $W_{\pm}$ depends on $\delta'$.) 
By inclusion-exclusion we can detect the condition that some $t_i  \in S_{\rho}(N)$ by an alternating sum of terms where a certain subset of the $t_i$ is in $S_{\rho}(N)$ while the complement is unrestricted. All of these terms are handled in the same way, so for notational simplicity let us focus on the case $t_1 \in S_{\rho}(N)$, $t_2, t_3, t_4$ unrestricted. We recall the notation $D' = D(1 + \delta')$ and define  
\begin{equation}\label{g}
G(y, D) = \sum_{ D < t \leq   D'} \frac{1}{t^{2 \pi i y}}, \quad G_{\rho, N}(y, D) = \sum_{\substack{ D < t \leq  D'\\ t \in S_{\rho}(N)}} \frac{1}{t^{2 \pi i y}}.
\end{equation}
 By partial summation we have
$$G(y, D) =  G^{\ast}(1 - 2\pi i y, D', D) + O(Y), \quad G^{\ast}(z, D_1, D_2) = \frac{D_1^{z} -D_2^{z}}{z} $$
and using \eqref{wein}, we also have
$$G_{\rho, N}(y, D) = (1-P_N) G^{\ast}(1 - 2\pi i y, D(1 + \eta'), D)  + O\big( YD \exp(- (\log N)^{1-2\rho})\big)$$
for $y \ll Y$. Let
$$\Phi_{\rho, N}(y, D_1, D_2) =  G_{\rho, N}(-y, D_1) G(-y, D_2), \quad \Phi(y, D_3, D_4) = G(y, D_3) G(y, D_4). $$
Detecting the congruence condition modulo 2, we have 
\begin{equation}\label{detect}
\begin{split}
& \sum_{\substack{t_i \in (D_i, D'_i]\\  t_1  \in S_{\rho}(N)\\ 2 \mid t_1 - t_2, t_3 - t_4}}  \Big( \frac{ t_1t_2}{t_3t_4}\Big)^{2 \pi i y} \\
  = &  \Big(\Phi_{\rho, X}(y, D_1, D_2) -\frac{\Phi_{\rho, X}(y, \frac{1}{2}D_1, D_2)  + \Phi_{\rho, X}(y, D_1, \frac{1}{2}D_2)}{2^{-2\pi i y}} +\frac{2 \Phi_{\rho, X}(y, \frac{1}{2}D_1, \frac{1}{2}D_2)}{4^{-2\pi i y}}\Big)\\
&   \Big(\Phi(y, D_3, D_4) -\frac{\Phi(y, \frac{1}{2}D_3, D_4)  + \Phi(y, D_3, \frac{1}{2}D_4)}{2^{2\pi i y}} +\frac{2 \Phi(y, \frac{1}{2}D_3, \frac{1}{2}D_4)}{4^{2\pi i y}}\Big) .
\end{split}
\end{equation}
Substituting all of this, we recast \eqref{eval} as
\begin{equation}\label{eval1}
\begin{split}
\frac{1}{2} &\sum_{(D_1, \ldots, D_4) \in \mathcal{D}_{+}} t\frac{4\sqrt{\alpha}\Delta}{\pi D_3D_4}  \widehat{W}_{+}(0)(1-P_N)\\ &\quad \int_{\Bbb{R}} V\Big(\frac{y}{Y}\Big) \frac{\alpha^{2\pi i y}}{2}
 \prod_{j=1}^2 G^{\ast}(1 + 2\pi i y, D_j', D_j)
  \prod_{j=3}^4 G^{\ast}(1 - 2\pi i y, D_j', D_j)
 dy \\
 &\quad + O\big(Y^2 \Delta N \exp(-(\log N)^{\rho/2})\big). 
\end{split}
\end{equation}
 Define 
$$H(t, D_1, D_2) = \delta_{\log D_2 < t \leq \log D_1} e^t.$$
Then $$\widehat{H}(y, D_1, D_2) = \int_{\Bbb{R}} H(t, D_1, D_2) e(ty) dt = G^{\ast}(1 + 2\pi i y, D_1, D_2).$$
Thus the main term in \eqref{eval1} becomes
\begin{displaymath}
\begin{split}
 \sum_{(D_1, \ldots, D_4) \in \mathcal{D}_{+}} \frac{\sqrt{\alpha}\Delta}{\pi D_3D_4}&  \widehat{W}_{+}(0)(1-P_N) \int_{\Bbb{R}} V\Big(\frac{y}{Y}\Big)  \int_{\log  D_1}^{\log  D_1'} \int_{\log  D_2}^{\log D_2'} \int_{\log D_3}^{\log D_3'}\\
& \int_{\log  D_4}^{\log  D_4'} e^{t_1+t_2+t_3+t_4}  e((\log \alpha + t_1 + t_2 - t_3 - t_4)y) dt_1\, dt_2\, dt_3\, dt_4\, dy . 
\end{split}
\end{displaymath}
Changing variables and recalling \eqref{defY}, this equals
\begin{equation}\label{thisequals}
\begin{split}
 \sum_{(D_1, \ldots, D_4) \in \mathcal{D}_{+}} & \frac{\sqrt{\alpha} \Delta Y'}{\pi N}  \widehat{W}_{+}(0)(1-P_N)   \int_{ D_1}^{  D_1'} \int_{  D_2}^{  D_2'} \int_{ D_3}^{ D_3'}  \int_{  D_4}^{  D_4'}   \widehat{V} \Big(Y \log\frac{ \alpha  t_1  t_2 }{t_3  t_4 }\Big) dt_1\, dt_2\, dt_3\, dt_4. 
\end{split}
\end{equation}
By the rapid decay of $\widehat{V}$ we can restrict to $\log \alpha  t_1  t_2 /t_3  t_4  \ll Y^{-3/4}$ at the cost of a total error $\Delta X Y^{-100}$. 
By a Taylor argument we then have
\begin{displaymath}
\begin{split}
&Y \log \frac{\alpha t_1t_2}{t_3t_4} = \frac{Y}{t_3t_4}(\alpha t_1 t_2 - t_3t_4) + O(Y^{-1/2})\\
& =\frac{ Y'}{ N} \frac{D_3D_4}{t_3t_4}(\alpha t_1 t_2 - t_3t_4) + O(Y^{-1/2})  = (1 + O(\delta')) \frac{ Y'}{ N} (\alpha t_1 t_2 - t_3t_4) + O(Y^{-1/2}).
\end{split}
\end{displaymath}
Using also $\widehat{W}_{+}(0)(1-P_X) = 1 + O(\delta')$  and defining 
$$\mathscr{D}_{+} = \bigcup_{(D_1, \ldots, D_4) \in \mathcal{D}_{+}} \bigtimes_{i=1}^4 [D_i, D_i'],$$
 we recast the main term \eqref{thisequals} as
\begin{displaymath}
\begin{split}
& \frac{(1 + O(\delta'))\sqrt{\alpha} \Delta Y' }{\pi N }    \int_{  \mathscr{D}_{+}}  \widehat{V} \Big(\frac{Y'}{N}  ( \alpha  t_1  t_2 - t_3  t_4 ) \Big)d\textbf{t} . 
\end{split}
\end{displaymath}
We now add to the  integration domain $\mathscr{D}_{+}$ those points of \eqref{region} with $\min(t_1, \ldots, t_4) \leq N^{1/2}(\log N)^{-A}$. It is easy to see (e.g.\ by putting the variables into dyadic boxes) that this infers an error of at most $O(\Delta N (\log N)^{2-A})$. Next we replace the integration domain with the exact region \eqref{region}, the error of which can be absorbed in the existing $(1 + O(\delta'))$-term. By symmetry we can assume $t_2 \geq t_1$ and $t_4 \geq t_3$ after multiplying by 4, and we drop the condition $t_1 \geq 0$ at the cost of dividing by 2 and note that for negative $t_1$ we automatically have $t_3 \leq 0$  up to a negligible error. 
Finally we reverse the change of variables \eqref{var} (the Jacobian infers a factor 2) getting
\begin{displaymath}
\begin{split}
& \frac{4 \sqrt{\alpha}}{\pi} (1 + O(\delta')) \Delta  \frac{Y'}{N}  \int \widehat{V} \Big(\frac{Y'}{N}( \alpha n_1^2 + m_1^2 - (\alpha n_2^2 +  m_2^2) ) \Big)dn_1\, dn_2\, dm_1\, dm_2 .
\end{split}
\end{displaymath}
where the integration is taken over
$$\alpha n_j^2 + m_j^2 \leq \frac{4\sqrt{\alpha} N}{\pi}, \quad n_1, n_2, m_1, m_2 \geq 0.$$ 
Changing variables, this equals
\begin{displaymath}
\begin{split}
&  (1 + O(\delta')) \Delta  \frac{Y'}{N}  \int_{0}^{N} \int_0^{N} \widehat{V} \Big(\frac{Y'}{N}(r_1 - r_2)\Big)dr_2\, dr_1\\
& = (1 + O(\delta')) \Delta  \frac{Y'}{N}  \int_{0}^{X} \int_{-r_1}^{N-r_1} \widehat{V} \Big(-\frac{Y'}{N}r_2\Big)dr_2\, dr_1.
\end{split}
\end{displaymath}
The portion $r_1 \leq N (Y')^{-3/4}$ and $r_1\geq N - N (Y')^{-3/4}$ is negligible, in the remaining part we can extend the $r_2$-integral to all of $\Bbb{R}$ by the rapid decay of $\widehat{V}$ and finally obtain
$$(1 + O(\delta')) \Delta \int_{N (Y')^{-3/4}}^{N - N (Y')^{-3/4}} V(0) dr_1 = ( 1+ O(\delta')) \Delta N$$
as desired.

  \section{Proof of Proposition \ref{prop6}}\label{proof6}
  
We have
\begin{displaymath}
\begin{split}
&\int_{\mathcal{R}_0} \Big|\frac{1}{2} \sum_{(D_1, \ldots, D_4) \in \mathcal{D}_{+}} \frac{4\sqrt{\alpha}\Delta}{\pi D_3D_4} I_2(\alpha)  \Big|^2 \frac{d\alpha}{\alpha}   \\
&\ll  \sup_{D_i}  \frac{\Delta^2 (\log\log N)^8}{D_3^2D_4^2}  \int_{\mathcal{R}_0}\Big( \int_{|y| \geq Y} \Big|\widehat{W}_{+}\Big(\frac{4\sqrt{\alpha}\Delta y}{\pi D_3D_4}\Big)\Big|\Big| \sum_{\substack{t_i \in (D_i, D'_i]\\ \text{some } t_i  \in S_{\rho}(N)\\ 2 \mid t_1 - t_2, t_3 - t_4}} \Big( \frac{\alpha t_1t_2}{t_3t_4}\Big)^{2 \pi i y}\Big| dy \Big)^2   \frac{d\alpha}{\alpha}. 
   \end{split}
\end{displaymath}
Again we treat the case $t_1 \in S_{\rho}(N)$ and $t_2, t_3, t_4$ unrestricted, the other cases being similar, and we note that $t_4$ is unrestricted in all cases.  We can detect the parity conditions in the same way as in \eqref{detect} which replaces potentially some $D_i$ by $\frac{1}{2}D_i$. Finally 
we replace the function  $\widehat{W}_{+}\big( 4\sqrt{\alpha}\Delta y/\pi D_3D_4\big)$ by a non-negative majorant $W^{\ast}(\Delta y/D_3D_4)$ that is independent of $\alpha$, where $W^{\ast}$ is some suitable Schwartz class function. Opening the square and integrating over $\alpha$, it suffices to bound 
\begin{equation}\label{bound}
\sup_{D_i}  \frac{\Delta^2 (\log N)^{\varepsilon} }{D_3^2D_4^2}    \int_{|y| \geq Y} \Big| W^{\ast}\Big(\frac{ \Delta y}{  D_3D_4}\Big)\Big|^2 |G_{\rho, N}(y, D_1) G(y, D_2) G(y, D_3)G(y, D_4)|^2
 dy 
\end{equation}
using the notation \eqref{g}. We now manipulate $G_{\rho, N}(y, D_1) $ similarly as in \cite[Lemma 12]{MR}. Let $I = (\exp((\log N)^{\rho})), \exp((\log N)^{1-\rho})]$ and $\omega_I(m)$ be the number of prime divisors of $m$ in $I$. We have
\begin{displaymath}
\begin{split}
G_{\rho, N}(y, D) & = \sum_{\substack{D < n \leq D' \\ n \in S_{\rho}(N)}} \frac{1}{n^{2\pi i y}} = \sum_{p\in I} \sum_{\substack{D/p < m \leq D'/p \\ m \in S_{\rho}(N)}} \frac{(\omega_I(m) + \delta_{(p,m) = 1})^{-1}}{(pm)^{2\pi i y}} \\ & = \sum_{p \in I} \sum_{\substack{D / p < m \leq D' / p \\ m \in S_{\rho}(N)}} \frac{(\omega_I(m) + 1)^{-1}}{(p m)^{2 \pi i y}} \\ & \qquad \qquad + \sum_{p \in I} \sum_{\substack{D / p < m \leq D' / p \\ m \in S_{\rho}(N)  ,  p \mid m}} \frac{\omega_I(m)^{-1}}{(p m)^{2\pi i y}} - \sum_{p \in I} \sum_{\substack{D / p < m \leq D' / p \\ m \in S_{\rho}(N),  p \mid m}} \frac{(\omega_{I}(m) + 1)^{-1}}{(p m)^{2\pi i y}} 
\end{split}
\end{displaymath}
(this is often called Ramar\'e's identity). We split the first $p$-sum into $O(\kappa^{-1} \log N)$ intervals of the shape $P< p \leq P(1+\kappa)$ for 
$$\kappa = (\log N)^{-C}$$
 and $\exp((\log N)^{\rho})) \leq P  \leq \exp((\log N)^{1-\rho})$. We argue as in  \cite[Lemma 12]{MR}  and write 
\begin{displaymath}
  \begin{split}
    & \sum_{p \in I} \sum_{\substack{D / p < m \leq D' / p \\ m \in S_{\rho}(N)}} \frac{(\omega_{I}(m) + 1)^{-1}}{(p m)^{2\pi i y}} \\ 
&= \sum_P \sum_{\substack{P < p \leq P(1 + \kappa)\\ p \in I}} \frac{1}{p^{2\pi i y}}  \sum_{\substack{D/(P(1+\kappa))< m \leq D'/P \\ D  \leq mp \leq D' \\ m \in S_{\rho}(N)}} \frac{(\omega_I(m) + 1)^{-1}}{m^{2\pi i y}} \\
& = \sum_P \sum_{\substack{P < p \leq P(1 + \kappa)\\ p\in I}} \frac{1}{p^{2\pi i y}}  \sum_{\substack{D/P < m \leq D'/P \\ m \in S_{\rho}(N)}} \frac{(\omega_I(m) + 1)^{-1}}{m^{2\pi i y}} + \sum_{m \in J} \frac{d_m}{m^{2\pi i y}} 
\end{split}
\end{displaymath}
for certain $|d_m| \leq 1$, where  $J = [D/(1+\kappa) , D(1+\kappa)] \cup [D', D'(1+\kappa)]$ 
and $P$ runs over a sequence of  the type $P_0(1+\kappa)^j$.   We substitute this back into \eqref{bound} getting
\begin{equation}\label{bound-new}
\begin{split}
& \sup_{D_i} \frac{\Delta^2(\log N)^{\varepsilon}  \kappa^{-1} \log N}{D_3^2D_4^2}  \sum_P  \int_{|y| \geq Y} \Big|W^{\ast}\Big(\frac{ \Delta y}{  D_3D_4}\Big)\Big|^2\\
 & \quad\quad\quad  |Q_P(y) R_P(y, D_1)  G(y, D_2) G(y, D_3)G(y, D_4)|^2 dy \\
  &+ \sup_{D_i}   \frac{\Delta^2 (\log N)^{\varepsilon}    }{D_3^2D_4^2}   \int_{\Bbb{R}} \Big|W^{\ast}\Big(\frac{ \Delta y}{  D_3D_4}\Big)  V(y, D_1) G(y, D_2)  G(y, D_3)G(y, D_4)\Big|^2 dy \\
 &+ \sup_{D_i}   \frac{\Delta^2 (\log N)^{\varepsilon}    }{D_3^2D_4^2}   \int_{\Bbb{R}} \Big|W^{\ast}\Big(\frac{ \Delta y}{  D_3D_4}\Big)  U(y, D_1) G(y, D_2)  G(y, D_3)G(y, D_4)\Big|^2 dy 
\end{split}
\end{equation}
where
\begin{displaymath}
\begin{split}
& Q_P(y)=  \sum_{\substack{P < p \leq P(1 + \kappa)\\ p\in I}} \frac{1}{p^{2\pi i y}} , \quad  R_P(y, D) =  \sum_{D/P < m \leq D'/P  } \frac{(\omega_I(m) + 1)^{-1}}{m^{2\pi i y}},\\
  & V(y, D) =  \sum_{m \in J } \frac{d_m}{m^{2\pi i y}} \\
  & U(y, D) = \sum_{p \in I} \sum_{\substack{D / p < m \leq D' / p \\ m \in S_{\rho}(N)  ,  p \mid m}} \frac{\omega_I(m)^{-1}}{(p m)^{2\pi i y}} - \sum_{p \in I} \sum_{\substack{D / p < m \leq D' / p \\ m \in S_{\rho}(N)  , p \mid m}} \frac{(\omega_{I}(m) + 1)^{-1}}{(p m)^{2\pi i y}} . 
\end{split}
\end{displaymath}
We estimate the three terms in \eqref{bound-new} separately and recall 
the  standard mean value estimate \cite[Theorem 9.1]{IK} for Dirichlet polynomials 
\begin{equation}\label{mean}
 \int_{-T}^T \Big|\sum_{n\leq X} \frac{a_n}{n^{it}}\Big|^2 dt \ll (T + X) \sum_{n \leq X} |a_n|^2.
\end{equation}
For  the second term in \eqref{bound-new} we apply \eqref{mean} directly getting the bound
\begin{equation}\label{bound0}
\begin{split}
& \sup_{D_i}   \frac{\Delta^2 (\log N )^{\varepsilon}   }{D_3^2D_4^2} \Big(\frac{D_3D_4}{\Delta} +  D_1D_2D_3D_4  \Big)  \sum_{m_1 \in J}\sum_{D_2 < m_2 \leq D_2'}\sum_{D_3 < m_2 \leq D_3'}\sum_{D_4 < m_2 \leq D_4'}\tau_4(m_1m_2m_3m_4)\\
 &\ll\sup_{D_i}   \frac{\Delta^2 (\log N   )^{\varepsilon} }{D_3^2D_4^2} \Big(\frac{D_3D_4}{\Delta} +  D_1D_2D_3D_4  \Big)   (D_1 D_2D_3D_4\kappa)^{1/2} \Big(\sum_{n \ll D_1D_2D_3D_4} \tau_4(n)\tau_4(n)^2\Big)^{1/2}\\
 & \ll  (\Delta N  + \Delta^2 N^2) (\log N)^{32 - C/2} \ll  \Delta^2 N^2 (\log N)^{32 - C/2}
\end{split}
\end{equation}
where we used Cauchy-Schwarz,  $\sum_{n \leq X} \tau_4(n)^3 \ll X (\log X)^{63} $ and the assumption $\Delta \geq N^{-1+\eta}$. 

For the third term in \eqref{bound-new} we proceed similarly, applying \eqref{mean} directly. This gives the bound
\begin{equation}
\begin{split}
& \sup_{D_i}   \frac{\Delta^2 (\log N )^{\varepsilon}   }{D_3^2D_4^2} \Big(\frac{D_3D_4}{\Delta} +  D_1D_2D_3D_4  \Big) \sum_{\substack{ p \in I \\ D_1 < m p^2 < D'_1}} \sum_{\substack{D_2 < m_2 \leq D_2' \\ D_3 < m_2 \leq D_3' \\ D_4 < m_2 \leq D_4' }}\tau_4(m_1 p^2 m_2m_3m_4)\\
 &\ll\sup_{D_i}   \frac{\Delta^2 (\log N   )^{\varepsilon} }{D_3^2D_4^2} \Big(\frac{D_3D_4}{\Delta} +  D_1D_2D_3D_4  \Big)  (D_1 D_2 D_3 D_4) \exp(- (\log N)^{\rho / 2}) \\
 & \ll  (\Delta N  + \Delta^2 N^2) \exp( - (\log N)^{\rho / 4})  \ll  \Delta^2 N^2 \exp( - (\log N)^{\rho / 4}) .
\end{split}
\end{equation}

In the first term in \eqref{bound-new} we split the $y$-integral into two parts. For those $y$ with $Q_P(y) \leq  P(\log N)^{-D}$ we use \eqref{mean}  to  bound their contribution by
\begin{equation}\label{bound1}
\begin{split}
& \sup_{D_i} \frac{\Delta^2 (\log N)^{\varepsilon} (\kappa^{-1} \log N)^2 P^2}{D_3^2D_4^2(\log N)^{2D} }  \Big( \frac{D_3D_4}{\Delta} + \frac{D_1D_2D_3D_4}{P}\Big) \sum_{n \ll D_1D_2D_3D_4/P} \tau_4(n)^2\\
 &\ll  (\Delta N P + \Delta^2 N^2) (\log N)^{2C + 18 - 2D} \ll \Delta^2 N^2 (\log N)^{2C + 18 - 2D}
 \end{split}
\end{equation}
since $\Delta \gg N^{-1+\eta}$ and $P \leq N^{\eta / 10}$. 

Now we treat the integral over the remaining $y$ where $Q_P(y) \geq  P(\log N)^{-D}$. Here we can discretize the integral and estimate it by a sum over certain points of distance at most 1. From \cite[Lemma 8]{MR} with $T = N^2$, $V = (\log N)^D$ and $\exp((\log N)^{\rho}) \leq P \leq \exp((\log N)^{1-\rho})$ we conclude the that number of such points is at most $\ll \exp((\log N)^{1-\rho/2})$. This gives the bound
\begin{equation}\label{final}
 \sup_{D_i} \frac{\Delta^2 (\log N)^{\varepsilon}  \kappa^{-1} \log N}{D_3^2D_4^2}  \sum_P \sum_{j} |Q_P(y_j) R_P(y_j, D_1)  G(y_j, D_2) G(y_j, D_3)G(y_j, D_4)|^2 
 \end{equation}
where $P$ runs over $O(\kappa^{-1} \log N)$ numbers and $j$ runs over $O(N^{\varepsilon})$ numbers with $Y \leq |y_j| \ll N^2$. We recall from \eqref{defY} that $Y \gg (\log N)^{B-2A}$. A standard application of Perron's formula and the convexity bound for the Riemann zeta function shows
$$G(y, D) \ll \frac{D \log N}{T}(T + |y|)^{\varepsilon} +  \frac{D}{1+|y|} + D^{1/2}(T + |y|)^{1/4+ \varepsilon}$$
for any parameter $T > 1$, so that under the present conditions $G(y_j, D_4) \ll D_4 Y^{-1}$. For the remaining sum over $j$ we use the discrete mean value theorem  \cite[Lemma 9]{MR} (which is \cite[Theorem 9.6]{IK}) and bound \eqref{final} by
\begin{equation}\label{bound2}
\begin{split}
& \sup_{D_i} \frac{\Delta^2 (\log N)^{\varepsilon}  (\kappa^{-1} \log N)^2}{D_3^2D_4^2}  \frac{D_4^2}{Y^2} (D_1D_2D_3 + \sqrt{D_1 D_2} N^{\varepsilon} ) \sum_{n \ll D_1D_2D_3} \tau_4(n)^2 \\
& \ll \Delta^2  N^2(\log N)^{18 - 2B + 4A + 2C}.
\end{split}
\end{equation}
Combining \eqref{bound0} --  \eqref{bound1} and \eqref{bound2} and choosing $A = 20$, $B = 400$, $C = 200$, $D = 300$, we complete the proof.

\end{document}